\documentstyle[amsfonts,12pt,myart]{article}

\normalbaselines
\font\tenbm=cmmib10
\font\sevenbm=cmmib7
\font\fivebm=cmmib5
\newfam\bmfam
\textfont\bmfam=\tenbm \scriptfont\bmfam=\sevenbm
\scriptscriptfont\bmfam=\fivebm
\newtheorem{theorem}{Theorem}[section]
\newtheorem{opred}{Definition}[section]
\begin{document}

\author{Yuri A. Rylov}
\title{Classification of Finite Subspaces of Metric Space Instead of Constraints on
Metric}
\date{Institute for Problems in Mechanics, Russian Academy of Sciences, 101-1,
Vernadskii Ave., Moscow, 117526, Russia\\
e-mail: rylov@ipmnet.ru}
\maketitle

\begin{abstract}
A new method of metric space investigation, based on classification of its
finite subspaces, is suggested. It admits to derive information on metric
space properties which is encoded in metric. The method describes geometry
in terms of only metric. It admits to remove constraints imposed usually on
metric (the triangle axiom and nonnegativity of the squared metric), and to
use the metric space for description of the space-time and other geometries
with indefinite metric. Describing space-time and using this method, one can
explain quantum effects as geometric effects, i.e. as space-time properties.
\end{abstract}

\section{Introduction}

Let $M=\{\rho ,\Omega\}$ be a metric space, where $\Omega$ is a set of
points, and $\rho $ be the metric, i.e.
\begin{equation}
\rho :\quad \Omega \times \Omega \rightarrow [0,\infty )\subset {\Bbb R}
\label{a1.1}
\end{equation}
\begin{equation}
\rho (P,P)=0,\qquad \rho (P,Q)=\rho (Q,P),\qquad \forall P,Q\in \Omega
\label{a1.2}
\end{equation}
\begin{equation}
\rho (P,Q)\geq 0,\qquad \rho (P,Q)=0,\quad {\rm \; \; \; \;} P=Q,\qquad
\forall P,Q\in \Omega  \label{a1.3}
\end{equation}
\begin{equation}
\rho (P,Q)+\rho (Q,R)\geq \rho (P,R),\qquad \forall P,Q,R\in \Omega
\label{a1.4}
\end{equation}

\begin{opred}
\label{d1} Any subset $\Omega ^{\prime }\subset \Omega $ of points of the
metric space $M=\{\rho ,\Omega \}$, equipped with the metric $\rho ^{\prime }
$ which is a contraction $\rho |_{\Omega ^{\prime }\times \Omega ^{\prime }}$
of the mapping (\ref{a1.1}). on the set $\Omega ^{\prime }\times \Omega
^{\prime }$ is called the metric subspace $M^{\prime }=\{\rho ^{\prime
},\Omega ^{\prime }\}$ of the metric space $M=\{\rho ,\Omega \}$.
\end{opred}

It is easy to see that the metric subspace $M^{\prime }=\{\rho ^{\prime
},\Omega ^{\prime }\}$ is a metric space.

\begin{opred}
\label{d5} The metric space $M_n({\cal P}^n)=\{\rho ,{\cal P}^n\}$ is called
a finite one, if it consists of a finite number of points ${\cal P}^n\equiv
\{P_i\}$, $i=0,1,\ldots n$.
\end{opred}

\begin{opred}
\label{d6} Finite metric space $M_n({\cal P}^n)=\{\rho ,{\cal P}^n\}$ is
called oriented one $\overrightarrow{M_n({\cal P}^n)}$, if the order of its
points ${\cal P}^n=\{P_0,P_1,\ldots P_n\}$ is given.
\end{opred}

\begin{opred}
\label{d7} An oriented finite subspace $\overrightarrow{M_n({\cal P}^n)}$ of
the metric space $M=\{\rho ,\Omega \}$ is called a multivector. It is
designed by means of $\overrightarrow{M_n({\cal P}^n)}\equiv \overrightarrow{%
{\cal P}^n}\equiv \overrightarrow{\{P_0,P_1,\ldots P_n\}}$
\end{opred}

\begin{opred}
\label{d2} A description is called $\sigma $-immanent one, if it does not
contain any references to objects or concepts other than subspaces of the
metric space or its metric.
\end{opred}

$\sigma$-immanency of a description means that it is immanent to the metric
space and it is carried out in terms of its metric and subspaces. The prefix
"$\sigma$" associates with the world function $\sigma$, which is connected
with the metric by means of the relation $\sigma =\rho ^2/2$. The name
"world function" was suggested by Singe \cite{S60}, who introduced it for
the Riemannian space description and used it for description of the
space-time geometry. A use of the world function $\sigma$ instead of the
metric $\rho$ appears to be more convenient.

The shortest, connecting two arbitrary points $P,Q\in\Omega$, is the basic
geometrical object which is constructed usually in the metric space $\{\rho
,\Omega\}$ \cite{ABN86}. One can construct an angle, triangle, different
polygons from segments of the shortest. Construction of two-dimensional and
three-dimensional planes in the metric space is rather problematic. At any
rate it is unclear how one could construct these planes, using the shortest
as the main geometrical object. It gives rise to think that the metric
geometry, i.e. geometry generated by the metric space, is less pithy, than
the geometry of the Euclidean space, where such geometric objects as
two-dimensional and three-dimensional planes can be build without any
problems. In fact it is not so. In the scope of the metric geometry one can
construct almost all geometric objects which can be constructed in the
Euclidean geometry, including $n$-dimensional planes. It is necessary to use
more effective method of the metric space description, than that based on
the use of the shortest.

Let ${\cal P}^n=\{P_0,P_1,\ldots ,P_n\}\subset\Omega$ be the set of $n+1$
points $P_i\in\Omega $, $(i=1,2,\ldots n)$ in $D$-dimensional Euclidean
space $\Omega ={\Bbb R}^D$, $(D>n)$. Let $\rho$ be the Euclidean metric in $%
{\Bbb R}^D$. Let us consider $(n+1)$-edr with vertices at the points ${\cal P%
}^n$. Its volume $S_n({\cal P}^n)$ may be presented in $\sigma$-immanent
form.
\begin{equation}
S_n({\cal P}^n)=\frac1{n!}\sqrt{F_n\left( {\cal P}^n\right)},  \label{a1.0}
\end{equation}
\begin{equation}
F_n\left( {\cal P}^n\right) =\det ||\left( {\bf P}_0{\bf P}_i.{\bf P}_0{\bf P%
}_k\right) ||,\qquad P_0,P_i,P_k\in \Omega ,\qquad i,k=1,2,...n  \label{a1.6}
\end{equation}
\begin{equation}
\left( {\bf P}_0{\bf P}_i.{\bf P}_0{\bf P}_k\right) \equiv \Gamma \left(
P_0,P_i,P_k\right) \equiv \sigma \left( P_0,P_i\right) +\sigma \left(
P_0,P_k\right) -\sigma \left( P_i,P_k\right) ,\qquad i,k=1,2,...n.
\label{a1.7}
\end{equation}
\begin{equation}
\sigma (P,Q)\equiv \frac 12\rho ^2(P,Q),\qquad \forall P,Q\in \Omega ,
\label{a1.8}
\end{equation}
where ${\bf P}_0{\bf P}_i$, $i=1,2,\ldots n$ are $n$ vectors of the
Euclidean space. If these $n$ vectors are linear independent, the volume $%
S_n({\cal P}^n)$ of $(n+1)$-edr does not vanish. If they are linear
dependent, $S_n({\cal P}^n)=0$. In virtue of (\ref{a1.0})-(\ref{a1.8}) the
condition $F_n\left( {\cal P}^n\right)=0$ is a $\sigma$-immanent criterion
of linear dependence of vectors ${\bf P}_0{\bf P}_i$, $i=1,2,\ldots n$.

The value of the $\sigma$-immanent function $F_n\left( {\cal P}^n\right)$ of
points ${\cal P}^n$ may serve as a criterion of "linear independence of
vectors" even in the case, when a linear vector space cannot be introduced
and the concept of linear independence cannot be defined via it. This shows
that the concept of linear independence is in reality something more
fundamental, than an attribute of a linear vector space. Essentially the
quantity $F_n\left( {\cal P}^n\right)$ is a characteristic of $(n+1)$-point
metric space $M_n=\{\rho ,{\cal P}^n\}$. This quantity appears, when one
identifies the quantity $\sqrt{F_n\left( {\cal P}^n\right)}$ with the length
$|M_n|\equiv |{\cal P}^n|$ of the finite metric space $M_n=\{\rho ,{\cal P}%
^n\}$.

Let $n$ vectors ${\bf P}_0{\bf P}_i$, $i=1,2,\ldots n$ are linear
independent. Then $F_n\left( {\cal P}^n\right)\ne 0$. Let us construct the
linear span of vectors ${\bf P}_0{\bf P}_i$, $i=1,2,\ldots n$, consisting of
vectors ${\bf P}_0{\bf R}$. Then $n+2$ vectors ${\bf P}_0{\bf P}_i$, $%
(i=1,2,\ldots n)$, ${\bf P}_0{\bf R}$ are linear dependent, and the point $R$
satisfies the $\sigma$-immanent relation $F_{n+1}\left( {\cal P}^n,
R\right)= 0$. The set ${\cal L}({\cal P}^n)=\{R|F_{n+1}\left( {\cal P}^n,
R\right)= 0\}$ of points $R$ is a $n$-dimensional plane, passing through $%
n+1 $ points ${\cal P}^n$. As far as this relation is $\sigma$-immanent, it
determines some set ${\cal T}({\cal P}^n)=\{R|F_{n+1}\left( {\cal P}^n,
R\right)= 0\} \subset\Omega$ of points in any metric space $M=\{\rho
,\Omega\}$. This set ${\cal T}({\cal P}^n)$, called the $n$th order tube, is
an analog of $n$-dimensional plane of the Euclidean space. $n+1$ points $%
{\cal P}^n$, determining the tube, will be referred to as basic points of
the tube, or its $(n+1)$-point $\sigma$-basis. The $n$th order tube may be
considered to be a natural geometric object (NGO) of the metric space,
determined by $(n+1)$-point metric space $\{\rho ,{\cal P}^n\}$. The $n$th
order tube is a $\sigma$-immanent geometric object.

Thus, always there exists a subspace of the metric space which is an analog
of $n$-dimensional Euclidean plane, but linear operations on vectors cannot
be defined always. This takes place, because the tubes ${\cal T}({\cal P}^n)$
have another structure than corresponding planes ${\cal L}({\cal P}^n)$ of
the Euclidean space. Let ${\cal Q}^n\subset {\cal L}({\cal P}^n)$ be other $%
(n+1)$-point $\sigma$-basis ($F_n\left( {\cal Q}^n\right)\ne 0$), belonging
to the plane ${\cal L}({\cal P}^n)$. Then ${\cal L}({\cal Q}^n)={\cal L}(%
{\cal P}^n)$. In the arbitrary metric space it is not so, and, in general. $%
{\cal T}({\cal Q}^n)\ne {\cal T}({\cal P}^n)$. In particular, in the
Euclidean space any two different points of a straight determine this
straight. In many cases for a metric space two basic points determine the
tube which coincides with the shortest (For instance, it is so for a
Riemannian space, considered as a metric space). Then any two points of the
shortest are the basic points of this shortest and determine it. But there
are cases, when it is not so. Then two different points other than basic
points determine a tube, but it is another tube. In the case of the second
order tubes (analog of two-dimensional Euclidean plane) the inequality $%
{\cal T}({\cal Q}^2)\ne {\cal T}({\cal P}^2)$ is more likely to be a rule
than an exception from the rule. For instance, it is true for the Riemannian
space considered as a metric one.

A possibility of the metric space description in terms of only the shortest
is restricted. Although exhibiting ingenuity, such a description may be
constructed. For instance, A.D.~Alexandrov showed that internal geometry of
two-dimensional boundaries of convex three-dimensional bodies may be
represented in the $\sigma$-immanent form \cite{A48}. Apparently, without
introducing tubes of the order higher than unity, the solution of similar
problem for three-dimensional boundaries of four-dimensional bodies is very
difficult.

The $\sigma$-immanent conception of the metric space description can be
formulated as a classification of all finite metric subspaces $M_n({\cal P}%
^n)=\{\rho ,{\cal P}^n\}$, $(n=1,2,\ldots )$, where ${\cal P}^n\subset\Omega$%
, is a set of $n+1$ points $P_i\in\Omega $, $(i=0,1,\ldots n)$. Any $M_n(%
{\cal P}^n)$ associates with the number $|M_n({\cal P}^n)|=|{\cal P}%
^n|=n!S_n({\cal P}^n)=\sqrt{F_n({\cal P}^n)}$, called the length.

From mathematical viewpoint the classification of metric subspaces reduces
to equipping the metric space $\{\rho ,\Omega\}$ with a series of $\sigma$%
-immanent mappings
\begin{equation}
F_n:\quad \Omega ^{n+1}\rightarrow {\Bbb R},\qquad \Omega
^{n+1}=\bigotimes\limits_{k=1}^{n+1}\Omega ,\qquad n=1,2,\ldots  \label{a1.5}
\end{equation}
where $F_n({\cal P}^n)$ is defined by the relations (\ref{a1.6})-(\ref{a1.8}%
). As one can see from (\ref{a1.6}), (\ref{a1.7}), in the case $n=1$ $%
F_1(P,Q)=2\sigma (P,Q)=\rho ^2(P,Q)$.

Further it will be shown that the classification of finite metric subspaces,
carried out by means of the series of mappings (\ref{a1.5}), admits to
derive information on the metric space properties contained in its metric.
The metric geometry (i.e. the geometry generated by the metric space)
appears to be not less pithy, than the Euclidean geometry. It means,
particularly, that the Euclidean geometry may be formulated in the $\sigma$%
-immanent form. Furthermore the geometry of any subset of points of the
proper Euclidean space may be formulated in the $\sigma$-immanent form. In
other words, the metric geometry, constructed on the basis of the metric, is
insensitive to continuity or discreteness of the space.

The situation of constructing the metric geometry may be presented
conveniently as follows. In the $D$-dimensional Euclidean space $E_D=\{\rho
,\Omega \}$, $\Omega ={\Bbb R}^D$ the $n$-dimensional plane ${\cal L}({\cal P%
}^{n})$, $n=1,2,\ldots D$, which passes through points ${\cal P}^n$, forming
$(n+1)$-point $\sigma$-basis $(F_n({\cal P}^n)\ne 0)$, is described as a set
of points $R$, satisfying $\sigma$-immanent equation $F_{n+1}({\cal P}^n,R)=
0$. $n$-dimensional plane ${\cal L}({\cal P}^{n})$, $(n=1,2,\ldots D)$ is
determined by only metric. It is NGO for the Euclidean space $E_D$.

The metric space can be conceived as a result of a deformation of $D$%
-dimensional Euclidean space $E_D$ with rather large $D$. The deformation
means a variation of distances between points of $E_D$, accompanied by
removing some set $U$ of points belonging to $E_D$. Under such a deformation
NGOs are deformed, turning to sets of points of more complicated
configuration, but they continue to be attributes of the metric space,
because they are $\sigma$-immanent and determined only by metric.

Restrictions (\ref{a1.3}), (\ref{a1.4}) on the metric $\rho$ are used by no
means. They are needed for constructing the shortest. They may be removed,
if geometrical objects are constructed on the basis of a classification of
finite metric subspaces $\{\rho ,{\cal P}^n\}$.

In this case, replacing the metric by the world function $\sigma =\frac
12\rho ^2$, one obtains the more general metric space $V=\{\sigma ,\Omega\}$
instead of the usual metric space $M=\{\rho ,\Omega\}$. This metric space
will be referred to as $\sigma $-space. The geometry, generated by the $%
\sigma $-space will be referred to as T-geometry. The T-geometry is a
generalization of the metric geometry on the case of indefinite metric.
T-geometry may be used for a description of the space-time geometry.

Under above described deformation of the Euclidean space the Euclidean
straights turn to hallow tubes. (In general, it is possible such a case,
when the straights turn to curves, remaining to be lines, but it ia a very
special case of deformation). The hallow tubes appear in the general case,
because one equation, determining the tube, describes generally a surface.
This fact explains the name of the geometry: tubular geometry, or
T-geometry. From viewpoint of the more general T-geometry the conventional
metric geometry is a degenerated geometry, where the tubes degenerate to
lines (the shortests). The T-geometry is a natural geometry (which is
nondegenerated, in general). It is the most general geometry. A strong
argument in favour of T-geometry is the circumstance that on the basis of
T-geometry one can construct such a space-time model, where quantum effects
are explained as simple T-geometric effects, and the quantum constant is an
attribute of the space-time \cite{R91}.

In the second section definitions of main objects of $\sigma $-space are
given. The third section is devoted to the formulation and proof of the
theorem, stating that the Euclidean geometry can be described in terms of
only metric. In the fourth section the role of the triangle axiom is
discussed.

\section{$\sigma $-space and its properties.}

\begin{opred}
\label{d3.1.1} $\sigma $-space $V=\{\sigma ,\Omega \}$ is nonempty set $%
\Omega $ of points $P$ with given on $\Omega \times \Omega $ real function $%
\sigma $
\begin{equation}
\sigma :\quad \Omega \times \Omega \to {\Bbb R},\qquad \sigma (P,P)=0,\qquad
\sigma (P,Q)=\sigma (Q,P)\qquad \forall P,Q\in \Omega .  \label{a2.1}
\end{equation}
\end{opred}

The function $\sigma $ is called world function, or $\sigma $-function.

\begin{opred}
\label{d3.1.1 }. Nonempty subset $\Omega ^{\prime }\subset \Omega $ of
points of the $\sigma $-space $V=\{\sigma ,\Omega \}$ with the world
function $\sigma ^{\prime }=\sigma |_{\Omega ^{\prime }\times \Omega
^{\prime }}$, which is a contraction $\sigma $ on $\Omega ^{\prime }\times
\Omega ^{\prime }$ is called $\sigma $-subspace $V^{\prime }=\{\sigma
^{\prime },\Omega ^{\prime }\}$ of $\sigma $-space $V=\{\sigma ,\Omega \}$.
\end{opred}

Further the world function $\sigma ^\prime = \sigma |_{\Omega ^{\prime
}\times\Omega ^\prime }$, which is a contraction of $\sigma $ will be
designed by means of $\sigma $. Any $\sigma$-subspace of $\sigma$-space is a
$\sigma$-space.

\begin{opred}
\label{d3.1.1ba}. $\sigma $-space $V^{\prime }=\{\sigma ^{\prime },\Omega
^{\prime }\}$ is called isometrically embedded in $\sigma $-space $%
V=\{\sigma ,\Omega \}$, if there exists such a monomorphism $f:\Omega
^{\prime }\rightarrow \Omega $, that $\sigma ^{\prime }(P,Q)=\sigma
(f(P),f(Q))$,\quad $\forall P,\forall Q\in \Omega ^{\prime },\quad
f(P),f(Q)\in \Omega $,
\end{opred}

Any $\sigma$-subspace $V^{\prime}$ of $\sigma$-space $V=\{\sigma ,\Omega \}$
is isometrically embedded in it.

\begin{opred}
\label{d3.1.1b}. Two $\sigma $-spaces $V=\{\sigma ,\Omega \}$ and $V^{\prime
}=\{\sigma ^{\prime },\Omega ^{\prime }\}$ are called to be isometric
(equivalent), if $V$ is isometrically embedded in $V^{\prime }$, and $%
V^{\prime }$ is isometrically embedded in $V$.
\end{opred}

\begin{opred}
\label{d3.1.1bc}. $\sigma $-space $M_n({\cal P}^n)=\{\sigma ,{\cal P}^n\}$,
consisting of $n+1$ points ${\cal P}^n$ is called the finite $\sigma $-space
of $n$th order.
\end{opred}

\begin{opred}
\label{d3.1.1bd}. The number $\sqrt{F_n({\cal P}^n)}$, where $F_n({\cal P}%
^n) $ is defined by relations (\ref{a1.6})-(\ref{a1.8}), is called the
length (volume) of the finite $n$th order $\sigma $-space $M_n({\cal P}^n)$.
\end{opred}

If the set of points ${\cal P}^n$ of a finite $\sigma$-space $M_n({\cal P}%
^n) $ is ordered, such a finite $\sigma$-space $M_n({\cal P}^n)$ is called
multivecor. Practically only multivectors which are $\sigma$-subspaces of
the same $\sigma$-space and described by the same world function will be
considered. In this case one may not mention on the world function in the
definition of the multivector and define it as follows.

\begin{opred}
\label{d3.1.2b}. The ordered set $\{P_l\},\quad l=0,1,\ldots n$ of $n+1$
points $P_0,P_1,...,P_n$, belonging to the $\sigma $-space $V$ is called the
$n$th order multivector $\overrightarrow{P_0P_1...P_n}$. The point $P_0$ is
the origin of the multivector $\overrightarrow{P_0P_1...P_n}$
\end{opred}

Let us use the following designation for the multivector $\overrightarrow{%
P_0P_1...P_n}$. $\overrightarrow{P_0P_1...P_n}\equiv {\bf P}_0{\bf P}_1...%
{\bf P}_n \equiv\overrightarrow{{\cal P}^n}$.

\begin{opred}
\label{d3.1.2}. The vector ${\bf PQ}$ in the $\sigma $-space $V$ is the
first order multivector, or the ordered set $\{P,Q\}$ of two points $P,Q$.
The point $P$ is the origin, nd $Q$ is the end of the vector.
\end{opred}

\begin{opred}
\label{d3.1.3}. The scalar $\sigma $-product $({\bf P}_0{\bf P}_1.{\bf P}_0%
{\bf P}_2)$ of two vectors ${\bf P}_0{\bf P}_1$ and ${\bf P}_0{\bf P}_2$,
having a common origin, is called a real number
\begin{equation}
({\bf P}_0{\bf P}_1.{\bf P}_0{\bf P}_2)\equiv \Gamma (P_0,P_1,P_2)\equiv
\sigma (P_0,P_1)+\sigma (P_0,P_2)-\sigma (P_1,P_2),  \label{a2.3}
\end{equation}
\[
P_0,P_1,P_2\in \Omega
\]
\end{opred}

In the case, when it does not lead to a misunderstanding, the term "scalar
product" will be used instead of the term "scalar $\sigma$-product".

\begin{opred}
\label{d3.1.4}. According to the definition \ref{d3.1.1bd}, the length $\mid
{\bf PQ}\mid $ of the vector ${\bf PQ}$ is the number
\begin{equation}
\mid {\bf PQ}\mid =\sqrt{2\sigma (P,Q)}=\left\{
\begin{array}{c}
\mid \sqrt{({\bf PQ.PQ})}\mid ,\quad ({\bf PQ.PQ})\ge 0 \\
i\mid \sqrt{({\bf PQ.PQ})}\mid ,\quad ({\bf PQ.PQ})<0
\end{array}
\right. \qquad P,Q\in \Omega  \label{a2.4}
\end{equation}
\end{opred}

\begin{opred}
\label{d3.1.5}. Vectors ${\bf P}_0{\bf P}_1$, ${\bf P}_0{\bf P}_2$ are
parallel or antiprallel, if the following relations are fulfilled
respectively
\begin{equation}
{\bf P}_0{\bf P}_1\uparrow \uparrow {\bf P}_0{\bf P}_2:\qquad ({\bf P}_0{\bf %
P}_1.{\bf P}_0{\bf P}_2)=\mid {\bf P}_0{\bf P}_1\mid \cdot \mid {\bf P}_0%
{\bf P}_2\mid  \label{a2.5}
\end{equation}
\begin{equation}
{\bf P}_0{\bf P}_1\uparrow \downarrow {\bf P}_0{\bf P}_2:\qquad ({\bf P}_0%
{\bf P}_1.{\bf P}_0{\bf P}_2)=-\mid {\bf P}_0{\bf P}_1\mid \cdot \mid {\bf P}%
_0{\bf P}_2\mid  \label{a2.6}
\end{equation}
\end{opred}

\begin{opred}
\label{d3.1.6}. Vectors ${\bf P}_0{\bf P}_1$, ${\bf P}_0{\bf P}_2$ are
collinear, if they are parallel, or antiparallel.
\begin{equation}
{\bf P}_0{\bf P}_1\parallel {\bf P}_0{\bf P}_2:\qquad ({\bf P}_0{\bf P}_1.%
{\bf P}_0{\bf P}_2)^2=\mid {\bf P}_0{\bf P}_1\mid ^2\cdot \mid {\bf P}_0{\bf %
P}_2\mid ^2  \label{a2.7}
\end{equation}
\end{opred}

\begin{opred}
\label{d3.1.6b}. The scalar $\sigma $-product $(\overrightarrow{{\cal P}^n}.%
\overrightarrow{{\cal Q}^n})$ of $n$th order multivectors $\overrightarrow{%
{\cal P}^n}$ and $\overrightarrow{{\cal Q}^n}$, having the common origin $%
P_0=Q_0$ is the real number
\begin{equation}
(\overrightarrow{{\cal P}^n}.\overrightarrow{{\cal Q}^n})=\det \Vert \Gamma
(P_0,P_i,Q_k)\Vert ,\,\,\,\,\,\,\,\,\,\,\,i,k=1,2,...n  \label{a2.8}
\end{equation}
\end{opred}

\begin{opred}
\label{d3.1.6c}. In accordance with the definition \ref{d3.1.1bd} the length
$|\overrightarrow{{\cal P}^n}|$ of the multivector $\overrightarrow{{\cal P}%
^n}$ is the number
\begin{equation}
|\overrightarrow{{\cal P}^n}|=\left\{
\begin{array}{c}
\mid \sqrt{(\overrightarrow{{\cal P}^n}.\overrightarrow{{\cal P}^n})}\mid =|%
\sqrt{F_n({\cal P}^n)}|,\quad (\overrightarrow{{\cal P}^n}.\overrightarrow{%
{\cal P}^n})\ge 0 \\
i\mid \sqrt{(\overrightarrow{{\cal P}^n}.\overrightarrow{{\cal P}^n})}\mid
=i|\sqrt{F_n({\cal P}^n)}|,\quad (\overrightarrow{{\cal P}^n}.%
\overrightarrow{{\cal P}^n})<0
\end{array}
\right. \qquad \overrightarrow{{\cal P}^n}\subset \Omega  \label{a2.9}
\end{equation}
where the quantity $F_n({\cal P}^n)$ is defined by the relatons (\ref{a1.6}%
)-(\ref{a1.8})
\end{opred}

\begin{opred}
\label{d3.1.5c}. Two $n$th order multivectors $\overrightarrow{{\cal P}^n}$ $%
\overrightarrow{{\cal Q}^n}$, having the common origin, are collinear $%
\overrightarrow{{\cal P}^n}\parallel \overrightarrow{{\cal Q}^n}$, if
\begin{equation}
(\overrightarrow{{\cal P}^n}.\overrightarrow{{\cal Q}^n})^2=|\overrightarrow{%
{\cal P}^n}|^2\cdot |\overrightarrow{{\cal Q}^n}|^2  \label{a2.10}
\end{equation}
\end{opred}

\begin{opred}
\label{d3.1.5e}. Two collinear $n$th order multivectors $\overrightarrow{%
{\cal P}^n}$ and $\overrightarrow{{\cal Q}^n}$ are similar oriented $%
\overrightarrow{{\cal P}^n}\uparrow \uparrow \overrightarrow{{\cal Q}^n}$
(parallel), if
\begin{equation}
(\overrightarrow{{\cal P}^n}.\overrightarrow{{\cal Q}^n})=|\overrightarrow{%
{\cal P}^n}|\cdot |\overrightarrow{{\cal Q}^n}|  \label{a2.11}
\end{equation}
They have opposite orientation $\overrightarrow{{\cal P}^n}\uparrow
\downarrow \overrightarrow{{\cal Q}^n}$ (antiparallel), if
\begin{equation}
(\overrightarrow{{\cal P}^n}.\overrightarrow{{\cal Q}^n})=-|\overrightarrow{%
{\cal P}^n}|\cdot |\overrightarrow{{\cal Q}^n}|  \label{a2.12}
\end{equation}
\end{opred}

{\it Example 2.1.} Let us consider $D$-dimensional point proper Euclidean
space. It may be considered as a metric space $E_D=\{\rho ,{\Bbb R}^D\}$ or
as a $\sigma $-space $E_D=\{\sigma ,{\Bbb R}^D\}$, the world function $%
\sigma =\frac 12\rho ^2$ being given by the relations
\begin{equation}
\sigma (P,Q)=\sigma (x,y)={\frac{1}{2}}\sum^{D}_{i,k=1}g_{ik}(x^{i}- y^{i})
(x^{k}- y^{k}), \qquad x,y\in {\Bbb R}^n,  \label{a2.13}
\end{equation}
\noindent where $x=\{x^{i}\}$ and $y=\{y^{i}\}$, $(i=1,2,\ldots D)$ are
contravariant coordinates of points $P$ and $Q$ respectively in some
rectilinear coordinate system $K$. Here $g_{ik}$ = const, $(i,k= 1,2,\ldots
D)$ is the metric tensor, $\det ||g_{ik}||\ne 0$. Eigenvalues of the matrix $%
g_{ik}$, $i,k= 1,2,\ldots D$ of the metric tensor are positive, and $%
\sum^{D}_{i,k=1}g_{ik}x^{i}x^{k}=0$, if and only if $x=0$. The above made
definitions of the vector, its length, scalar product of two vectors and
relations of collinearity agree with the use of these concepts for the
Euclidean space. Indeed, the length of the vector ${\bf PQ}$ in the
Euclidean space is

\begin{equation}
\mid {\bf PQ}\mid = \sqrt{g_{ik}(x^i-y^i)(x^k-y^k)}= \sqrt{2\sigma (P,Q)}
\label{a2.14}
\end{equation}
\noindent that agrees with (\ref{a2.4}).

In the proper Euclidean space according to the cosine theorem for two
vectors ${\bf P}_{0}{\bf P}_{1}$ and ${\bf P}_{0}{\bf P}_{2}$
\begin{equation}
\mid {\bf P}_{1}{\bf P}_{2}\mid ^{2}=\mid {\bf P}_{0}{\bf P}_{2}-{\bf P}_{0}%
{\bf P}_{1}\mid ^{2}=\mid {\bf P}_{0}{\bf P}_{2}\mid ^{2}+\mid {\bf P}_{0}%
{\bf Q}_{1}\mid ^{2}-2({\bf P}_{0}{\bf P}_{1}.{\bf P}_{0}{\bf P}_{2})
\label{a2.15}
\end{equation}
It follows from this relation
\begin{equation}
({\bf P}_{0}{\bf P}_{1}.{\bf P}_{0}{\bf P}_{2})={\frac{1}{2}}\{\mid {\bf P}%
_{0}{\bf P}_{2}\mid ^{2}+\mid {\bf P}_{0}{\bf P}_{1}\mid ^{2}-\mid {\bf P}%
_{1}{\bf P}_{2}\mid ^{2}\}  \label{a2.16}
\end{equation}
\noindent that agrees with (\ref{a2.3}), if one takes into account (\ref
{a2.4}).

In the proper Euclidean space the vectors ${\bf P}_{0}{\bf P}_{1}$ and ${\bf %
P}_{0}{\bf P}_{2}$ are parallel or antiparallel, if cosine of the angle $%
\vartheta $ between them is equal respectively to $1$ or $-1$. As far as
\begin{equation}
\cos \vartheta =({\bf P}_{0}{\bf P}_{1}.{\bf P}_{0}{\bf P}_{2}) \mid {\bf P}%
_{0}{\bf P}_{2}\mid ^{-1}\cdot\mid {\bf P}_{0}{\bf P}_{1} \mid ^{-1},
\label{a2.17}
\end{equation}
one obtains an accord with the definitions (\ref{a2.5}), (\ref{a2.6}).

In the proper Euclidean space the $n$th order multivector ${\bf m}$ is
defined as an external (skew) product of vectors
\begin{equation}
{\bf m=e}_1\wedge {\bf e}_2\wedge ...\wedge {\bf e}_n\qquad {\bf e}_i={\bf P}%
_0{\bf P}_i,\qquad i=1,2,...n  \label{a2.18}
\end{equation}
The scalar product of two $n$th order multivectors ${\bf m}$ and ${\bf q}$
\begin{equation}
{\bf q}=\bigwedge\limits_{i=1}^{i=n}{\bf k}_i,\qquad {\bf k}_i={\bf P}_0{\bf %
Q}_i{\bf \qquad }i=1,2,...n  \label{c.1.12b}
\end{equation}
is defined by means of the relation
\begin{equation}
\left( {\bf m.q}\right) =\det \left\| \left( {\bf e}_i.{\bf k}_l\right)
\right\| =\det \left\| \left( {\bf P}_0{\bf P}_i.{\bf P}_0{\bf Q}_l\right)
\right\| =\det \left\| \Gamma \left( P_0,P_i,Q_l\right) \right\| ,\qquad
i,l=1,2,...n,  \label{a2.19}
\end{equation}
that agrees with the relation (\ref{a2.8}).

The difference between the definition \ref{d3.1.2b} of the multivector and
its conventional definition (\ref{a2.18}) consists in that that the first
definition does not use the summation operation and that of multiplication
of vectors by a number which are not defined in $\sigma $-space, whereas the
conventional definition (\ref{a2.18}) refers to the concept of manifold and
linear space, where these operations are defined.

In principle, summation of vectors and multiplication of them by a number
may be defined in $\sigma$-space $V=\{\sigma ,\Omega \}$ as follows. The
vector ${\bf P}_0{\bf R}$ is a sum of vectors ${\bf P}_0{\bf P}_1$ and ${\bf %
P}_0{\bf P}_2$, if $\exists R\in\Omega $ such, that
\begin{equation}
\left( {\bf P}_0{\bf R}.{\bf P}_0{\bf Q}\right) =\left( {\bf P}_0{\bf P}_1.%
{\bf P}_0{\bf Q}\right) +\left( {\bf P}_0{\bf P}_2.{\bf P}_0{\bf Q}\right)
,\qquad \forall Q\in \Omega .  \label{a2.20}
\end{equation}
The vector ${\bf P}_0{\bf R}$ is a result of multiplication of ${\bf P}_0%
{\bf P}$ by the real number $a$:\quad ${\bf P}_0{\bf R}=a{\bf P}_0{\bf P}$,
if $\exists R\in \Omega $ such, that
\begin{equation}
\left( {\bf P}_0{\bf R}.{\bf P}_0{\bf Q}\right) =a\left( {\bf P}_0{\bf P}.%
{\bf P}_0{\bf Q}\right) ,\qquad \forall Q\in \Omega  \label{a2.21}
\end{equation}
But such definitions are not effective, because in general case there is no
point $R$, satisfying the relations (\ref{a2.20}), (\ref{a2.21}). In the
case of the proper Euclidean space, when the points $R$, satisfying (\ref
{a2.20}), (\ref{a2.21}), exist, the operation of summation of vectors (\ref
{a2.20}) and that of multiplication of the vector by a number (\ref{a2.21})
coincide with the conventional definition of these operations in the
Euclidean space.

Operation of permutation of the multivector points can be effctively defined
in the $\sigma $-space. Let us consider two $n$th order multivectors $%
\overrightarrow{{\cal P}^n}=\overrightarrow{P_0P_1P_2...P_n}$ and $%
\overrightarrow{{\cal P}_{(1\leftrightarrow 2)}^n}= \overrightarrow{%
P_0P_2P_1P_3P_4...P_n}$, $(n\ge 2)$, which differ by the order of points $%
P_1 $ ¨ $P_2$.
\begin{equation}
( \overrightarrow{{\cal P}^n}.\overrightarrow{{\cal Q}^n} ) =\det \parallel
\Gamma ( P_0,P_{i},Q_k)\parallel , \qquad i,k=1,2,...n, \qquad \forall {\cal %
Q}^n\subset\Omega ,  \label{a2.22}
\end{equation}
\begin{equation}
( \overrightarrow{{\cal P}_{(1\leftrightarrow 2)}^n}.\overrightarrow{{\cal Q}%
^n}) =\det \parallel \Gamma ( P_0,P_{i}^{\prime },Q_k) \parallel ,\qquad
i,k=1,2,...n, \qquad n\ge 2, \qquad \forall {\cal Q}^n\subset\Omega .
\label{a2.23}
\end{equation}
Here the ordered set of points $\{P_i^{\prime }\},\,\,\,\,\,\,\,
(i=1,2,...n) $ is obtained from the ordered set of points $%
\{P_i\},\,\,\,\,\,\,\,(i=1,2,...n)$ by permutation of points $P_1$ ¨ $P_2$.
This means that the determinant (\ref{a2.23}) is obtained from the
determinant (\ref{a2.22}) by permutation of the first and second rows. Then
one obtains
\begin{equation}
( \overrightarrow{{\cal P}^n}.\overrightarrow{{\cal Q}^n}) =-(
\overrightarrow{{\cal P}_{(1\leftrightarrow 2)}^n}.\overrightarrow{{\cal Q}^n%
}), \qquad n\ge 2,\qquad \forall {\cal Q}^n\subset\Omega .  \label{a2.24}
\end{equation}
As far as $\overrightarrow{{\cal Q}^n}$ is an arbitrary multivector, the
relation (\ref{a2.24}) may be written in the form
\begin{equation}
\overrightarrow{{\cal P}^n} =-\overrightarrow{{\cal P}_{(1\leftrightarrow
2)}^n}, \qquad n\ge 2.  \label{a2.25}
\end{equation}
It may be interpreted in the sense, that permutation of any two points $P_i$
and $P_k\,\,\,\,\,\,i,k=1,2,...n$, $(n \ge 2)$ (except for the origin $P_0$
) at the multivector $\overrightarrow{{\cal P}^n}$ leads to a change of its
sign. The negative sign of the multivector means by definition that
\begin{equation}
(- \overrightarrow{{\cal P}^n}.\overrightarrow{{\cal Q}^n}) =-(
\overrightarrow{{\cal P}^n}.\overrightarrow{{\cal Q}^n}) , \qquad \forall
{\cal Q}^n\subset\Omega .  \label{a2.26}
\end{equation}

Permutating the points $P_0$ and $P_1\,\,\,\,\,$ at the multivector $%
\overrightarrow{{\cal P}^n}$, $(n\ge 2)$, one turns it in multivector $%
\overrightarrow{{\cal P}_{(0\leftrightarrow 1)}^n}$, having the origin at
the point $P_1$. Strictly, one cannot compare multivectors $\overrightarrow{%
{\cal P}^n}$ and $\overrightarrow{{\cal P}_{(0\leftrightarrow 1)}^n}$ at the
point $P_0$. But they have other common points $P_2,P_3,...P_n$, and one may
compare them at these points, forming scalar product with the multivector $%
\overrightarrow{{\cal Q}^n}$, having the origin, for instance, at the point $%
P_2$ $(Q_2=P_2)$
\begin{equation}
( \overrightarrow{{\cal P}^n}.\overrightarrow{{\cal Q}^n}) _{P_2}=\det
\left\| \Gamma \left( P_2,P_i,Q_k\right) \right\|
,\,\,\,\,\,\,\,\,\,\,\,i,k=0,1,3,4,...n,\qquad n\ge 2,  \label{a2.27}
\end{equation}
\begin{equation}
( \overrightarrow{{\cal P}_{(0\leftrightarrow 1)}^n}.\overrightarrow{{\cal Q}%
^n}) _{P_2}=\det\parallel\Gamma ( P_2,P_{i}^{\prime }, Q_k)\parallel,\qquad
i,k=0,1,3,4,...n,\qquad n\ge 2,  \label{a2.28}
\end{equation}
where $P_0^{\prime }=P_1$,\,\, $P_1^{\prime }=P_0,$ $P_i^{\prime}=P_i,\,\,\,%
\,\,\,i=3,4,...n$, and index $P_2$ shows, that the point $P_2$ is considered
as the origin of the multivector $\overrightarrow{{\cal Q}^n}$. Comparison
of rhs of (\ref{a2.27}) and (\ref{a2.28}) shows that
\[
( \overrightarrow{{\cal P}^n}.\overrightarrow{{\cal Q}^n}) _{P_2}=-(
\overrightarrow{{\cal P}_{(0\leftrightarrow 1)}^n}.\overrightarrow{{\cal Q}^n%
}) _{P_2},\qquad n\ge 2.
\]
The same result is obained, choosing any point of $P_i\,\,\,\,\,\,\,\,%
\,i=3,4,...n$ as an origin. It means that the relation (\ref{a2.25}) is
valid for permutation of any two points of the multivector $\overrightarrow{%
{\cal P}^n}$, and one may write
\begin{equation}
\overrightarrow{{\cal P}_{(i\leftrightarrow k)}^n}=-\overrightarrow{{\cal P}%
^n},\qquad i,k=0,1,...n,\qquad i\neq k,\qquad n\ge 2.  \label{a2.29}
\end{equation}
Thus, a change of the $n$th order multivector sign $(n\ge 2)$
(multiplication by the number $a=-1$) may be always defined as an odd
permutation of points.

For the vector (the first order multivector) the multiplication (\ref{a2.21}%
) by the number $a=-1$ is defined directly as a permutation of the origin
and the end of the vector by means of the relatons
\[
-{\bf P}_0{\bf P}_1={\bf P}_1{\bf P}_0,
\]
It means by definition that
\[
\left( -{\bf P}_0{\bf P}_1.{\bf P}_0{\bf Q}\right) =-\left( {\bf P}_0{\bf P}%
_1.{\bf P}_0{\bf Q}\right)= -\sigma \left( P_0,P_1\right) -\sigma \left(
P_0,Q\right) +\sigma \left( P_1,Q\right) ,\qquad \forall Q\in \Omega ,
\]
\[
\left( -{\bf P}_0{\bf P}_1.{\bf P}_1{\bf Q}\right) =\left( {\bf P}_1{\bf P}%
_0.{\bf P}_1{\bf Q}\right)= \sigma \left( P_1,P_0\right) +\sigma \left(
P_1,Q\right) -\sigma \left( P_0,Q\right) ,\qquad \forall Q\in \Omega .
\]
Thus multiplication of any multivector by the number $a=\pm 1$ may be always
defined in $\sigma$-space as a result of permutation of points, forming the
multivector.

In the properEuclidean space, where the multivector is defined in the form
\begin{equation}
\overrightarrow{{\cal P}^n}=\bigwedge_{i=1}^{i=n}{\bf P}_0{\bf P}_i ,
\label{a2.30}
\end{equation}
it is antisymmetric with respect to permutation of any two indices $%
i,k=0,1,...n$,\quad $i\neq k$. For indices $i,k=1,2,...n,$ it follows from
the external product properties.

For permutation of points $P_0\leftrightarrow P_1$ one has
\[
\overrightarrow{{\cal P}_{(0\leftrightarrow 1)}^n}={\bf P}_1{\bf P}%
_0\bigwedge_{i=2}^{i=n}{\bf P}_1{\bf P}_i=-{\bf P}_0{\bf P}%
_1\bigwedge_{i=2}^{i=n}\left( {\bf P}_0{\bf P}_i-{\bf P}_0{\bf P}_1\right) =
\]
\begin{equation}
-{\bf P}_0{\bf P}_1\bigwedge_{i=2}^{i=n}{\bf P}_0{\bf P}_i=-\overrightarrow{%
{\cal P}^n}  \label{a2.31}
\end{equation}
A similar result is obtained for permutation of points $P_0\leftrightarrow
P_i,\,\,\,\,\,\,\,i=1,2,...n$. Thus, the multivector in $\sigma $-space is
the geometrical object antisymmetric with respect to permutation of any two
points.

\begin{opred}
\label{d3.1.7}. $n+1$ points ${\cal P}^n$ , $P_i\in \Omega \quad (i=0,1,..n)$
form $(n+1)$-point $\sigma $-basis of the tube in $\sigma $-space, if the
multivector $\overrightarrow{{\cal P}^n}$ has nonvanishing length
\begin{equation}
|\overrightarrow{{\cal P}^n}|^2\equiv F_n({\cal P}^n)\ne 0.  \label{a2.32}
\end{equation}
\end{opred}

Let us illustrate this definition of the tube $\sigma $-basis in the example
of the $D$-dimensional proper Euclidean space. Let $n$ vectors ${\bf e}_{i}=%
{\bf P}_{0}{\bf P}_{i}$, $i=1,2,\ldots n$ be given in $D$-dimensional proper
Euclidean space $(n\le D)$. In this case (\ref{a1.6}) is the Gram's
determinant
\begin{equation}
F_{n}({\cal P}^{n})=\det \parallel ({\bf e}_{i}.{\bf e}_{k})\parallel
=(n!S_{n}({\cal P}^{n}))^{2},\qquad i,k=1,2,\ldots n  \label{a2.33}
\end{equation}
\noindent and $S_{n}$ is the volume of $(n+1)$-edr with vertices at points $%
{\cal P}^{n}$. Vanishing of this determinant is the necessary and sufficient
condition of linear independence of vectors ${\bf e}_{i},\quad (i=1,2,\ldots
n)$ in the proper Euclidean space.

If the condition (\ref{a2.32}) is fulfilled, $n$ vectors ${\bf e}_{i}$ are
linear independent and may serve as a basis in the $n$-dimensional plane $%
{\cal L}({\cal P}^{n})$, passing through points ${\cal P}^{n}$. In
particular, if one uses the expression (\ref{a2.13}) for calculation of the
scalar product of the vectors ${\bf e}_{i}={\bf P}_{0}{\bf P}_{i}$,\quad $%
(i=1,2,\ldots D)$, considering the $(D+1)$-point tube $\sigma$-basis ${\cal P%
}^{D}$, as the system of coordinate vectors, one obtains by means of (\ref
{a2.3}) ¨ (\ref{a2.13}) that
\begin{equation}
({\bf e}_{i}.{\bf e}_{k})=({\bf P}_{0}{\bf P}_{i}.{\bf P}_{0}{\bf P}_{k})
=g_{ik}({\cal P}^{D})=\Gamma (P_{0},P_{i},P_{k}),\qquad i,k=1,2,\ldots D
\label{a2.34}
\end{equation}

\begin{opred}
\label{d3.1.8}. The $n$th order tube ${\cal T}({\cal P}^n)$, (n=0,1,\ldots
), formed by $(n+1)$-point tube $\sigma $-basis ${\cal P}^n\subset \Omega $
(or by the $n$th order multivector $\overrightarrow{{\cal P}^n}\subset
\Omega $), is the set of points $P\in \Omega $
\begin{equation}
{\cal T}({\cal P}^n)\equiv {\cal T}_{{\cal P}^n}=\{P\mid F_{n+1}(P,{\cal P}%
^n)=0\},\qquad F_n({\cal P}^n)\neq 0.  \label{a2.35}
\end{equation}
\end{opred}

The relation (\ref{a2.35}) may be written also in terms of multivector $%
\overrightarrow{{\cal P}^n}$
\begin{equation}
{\cal T}({\cal P}^{n}) = \left\{P_{n+1}\left| |\overrightarrow{{\cal P}^{n+1}%
} |=0\right.\right\}, \qquad | \overrightarrow{{\cal P}^{n}} |\neq 0.
\label{a2.36}
\end{equation}

The tube ${\cal T}({\cal P}^{n})$ is the $n$th order natural geometrical
object (NGO), i.e the set of points, determined by geometry and parameters: $%
n+1$ points ${\cal P}^{n}$. The set of all possible NGOs is a set of $\sigma$%
-immanent geometric objects on the set $\Omega$. Each NGO contains at least
basic points ${\cal P}^{n}$.

\begin{opred}
\label{d3.1.9}. Section ${\cal S}_{n;P}$ of the tube ${\cal T}({\cal P}^n)$
at the point $P\in {\cal T}({\cal P}^n)$ is the set ${\cal S}_{n;P}({\cal T}(%
{\cal P}^n))$ of points, belonging to the tube ${\cal T}({\cal P}^n)$
\begin{equation}
{\cal S}_{n;P}({\cal T}({\cal P}^n))=\{P^{\prime }\mid
\bigwedge_{l=0}^{l=n}\sigma (P_l,P^{\prime })=\sigma (P_l,P)\},\qquad
P,P^{\prime }\in {\cal T}({\cal P}^n).  \label{a2.38}
\end{equation}
\end{opred}

In the proper Euclidean space the $n$th order tube is the $n$-dimensional
plane, containing points ${\cal P}^{n}$ , and its sectiom ${\cal S}_{n;P}(%
{\cal T}({\cal P}^{n}))$ at the point $P$ consists of one point $P$.

The zeroth and first order tubes are the most interesting and important. For
$F_{1}$ one obtains from (\ref{a1.6}) and (\ref{a2.3})
\[
F_{1}(P_{0},P_{1}) = 2\sigma (P_{0},P_{1})
\]
Then
\begin{equation}
{\cal T}(P_{0})\equiv {\cal T}_{P_{0}}=\{P\mid\sigma (P_{0},P)=0\},
\label{a2.39}
\end{equation}
In the properEuclidean space the zeroth order tube ${\cal T}_{P_{0}}
=\{P_{0}\}$ consists of one point $P_{0}$, and its section ${\cal S}%
_{0;P_{0}}({\cal T}_{p_{0}}) =\{P_{0}\}$ consists of one point $P_{0}$ also.
But in the pseudo-Euclidean space (for instance, in the space-time of the
special relativity) ${\cal T}_{P_{0}}$ is the light cone with the vertex at
the point $P_{0}$, and its section
\[
{\cal S}_{0;P}({\cal T}(P_{0}))=\{P^{\prime}\mid \sigma
(P_{0},P^{\prime})=0\wedge \sigma (P_{0},P^{\prime})=\sigma (P_{0},P)\}=%
{\cal T}_{P_{0}}
\]
\noindent coincides with the light cone.

Describing the first order tubes, it is convenient to use the circumstance
that the function $F_{2}({\cal P}^{2})$ can be presented in the form of a
product
\begin{equation}
F_{2}(P_{0},P_{1},P_{2})=S_{+}(P_{0},P_{1},P_{2})S_{2}(P_{0},P_{1},P_{2})S_{2}(P_{1},P_{2},P_{0})S_{2}(P_{2},P_{0},P_{1})
\label{a2.40}
\end{equation}
\noindent where
\begin{equation}
S_{+}(P_{0},P_{1},P_{2})\equiv S(P_{0},P_{1})+S(P_{1},P_{2})+S(P_{0},P_{2})
\label{c.1.22}
\end{equation}
\begin{equation}
S_{2}(P_{0},P_{1},P_{2})\equiv S(P_{0},P_{1})+S(P_{1},P_{2})-S(P_{0},P_{2})
\label{c.1.23}
\end{equation}
Here $S=\sqrt{2\sigma }$. $S_{+}$ vanishes, if and only if any term of the
sum (\ref{c.1.22}) vanishes. Then no two points form $\sigma$-basis, and the
tube is not defined. The tube ${\cal T}({\cal P}^{2})$ may be presented as
consisting of parts, and any multiplier in (\ref{a2.40}) (except for $S_{+}$%
) is responsible for one of these parts.

Let us set
\begin{equation}
{\cal T}_{[P_{0}P_{1}]}={\cal T}_{[P_{1}P_{0}]}=\{P\mid
S_{2}(P_{0},P,P_{1})=0\}  \label{a2.41}
\end{equation}
\begin{equation}
{\cal T}_{P_{0}[P_{1}}={\cal T}_{P_{1}]P_{0}}=\{P\mid
S_{2}(P_{0},P_{1},P)=0\}  \label{c.1.25}
\end{equation}
Let us refer to ${\cal T}_{[P_{0}P_{1}]}$ as the tube segment between the
points $P_{0}$, $P_{1}$, and to ${\cal T}_{P_{0}[P_{1}}$ as the tube ray
outgoing from $P_{1}$ towards th point $P_{0}$.

It is evident from (\ref{a2.40}), (\ref{a2.41}), (\ref{c.1.25}) that
\begin{equation}
{\cal T}_{P_{0}P_{1}}={\cal T}_{P_{0}]P_{1}}\bigcup {\cal T}_{[P_{0}P_{1}]}
\bigcup {\cal T}_{P_{0}[P_{1}}  \label{a2.42}
\end{equation}
As far as the relation (\ref{a2.7}) is equivalent to the equation $F_{2}(%
{\cal P}^{2})=F_{2}(P_{2},{\cal P}^{1})=0$, the first order tube ${\cal T}%
_{P_{0}P_{1}}$ may be defined also as a set of such points $P$ that ${\bf P}%
_{0}{\bf P}\parallel {\bf P}_{0}{\bf P}_{1}$.
\begin{equation}
{\cal T}_{P_{0}P_{1}}=\{P\mid {\bf P}_{0}{\bf P}_{1}\parallel {\bf P}_{0}%
{\bf P\}}  \label{a2.43}
\end{equation}
For the tube rays one can use the definitions
\begin{equation}
{\cal T}_{P_{0}[P_{1}}=\{P\mid {\bf P}_{1}{\bf P}\uparrow\downarrow {\bf P}%
_{1}{\bf P}_{0}\}  \label{a2.44}
\end{equation}
\begin{equation}
{\cal T}_{[P_{0}P_{1}}=\{P\mid {\bf P}_{0}{\bf P}\uparrow\uparrow {\bf P}_{0}%
{\bf P}_{1}\}  \label{a2.45}
\end{equation}

\begin{opred}
\label{d3.1.8 }. The oriented segment $\overrightarrow{{\cal T}_{[P_0P_1]}}$
of the first order tube, formed by by the vector $\overrightarrow{P_0P_1}%
\subset \Omega $ of unvinishing length is a totality $\{\overrightarrow{%
P_0P_1},{\cal T}_{[P_0P_1]}\}$ of the vector $\overrightarrow{P_0P_1}$ and
segment ${\cal T}_{[P_0P_1]}$, formed by this vector. The length of the
oriented segment $\overrightarrow{{\cal T}_{[P_0P_1]}}$ is the quantity
\begin{equation}
|\overrightarrow{{\cal T}_{[P_0P_1]}}|=|\overrightarrow{P_0P_1}|=\sqrt{%
2\sigma (P_0,P_1)}.  \label{a2.46}
\end{equation}
The scalar $\sigma $-product $(\overrightarrow{{\cal T}_{[P_0P_1]}}.%
\overrightarrow{P_0Q})$ of the oriented segment $\overrightarrow{{\cal T}%
_{[P_0P_1]}}$ and vector $\overrightarrow{P_0Q}$ is the number
\begin{equation}
(\overrightarrow{{\cal T}_{[P_0P_1]}}.\overrightarrow{P_0Q})=(%
\overrightarrow{P_0P_1}.\overrightarrow{P_0Q})=\sigma (P_0,P_1)+\sigma
(P_0,Q)-\sigma (P_1,Q),\qquad P_0,P_1,P,Q\in \Omega .  \label{a2.47}
\end{equation}
The scalar $\sigma $-product $(\overrightarrow{{\cal T}_{[P_0P_1]}}.%
\overrightarrow{P_1Q})$ of the oriented segment $\overrightarrow{{\cal T}%
_{[P_0P_1]}}$ and vector $\overrightarrow{P_1Q}$ is the number
\begin{equation}
(\overrightarrow{{\cal T}_{[P_0P_1]}}.\overrightarrow{P_1Q})=-(%
\overrightarrow{P_1P_0}.\overrightarrow{P_1Q})=-\sigma (P_0,P_1)-\sigma
(P_1,Q)+\sigma (P_0,Q),\qquad P_0,P_1,P,Q\in \Omega .  \label{a2.48}
\end{equation}
\end{opred}

In other words, $\overrightarrow{{\cal T}_{[P_0P_1]}}=-\overrightarrow{{\cal %
T}_{[P_1P_0]}}$.

Describing in terms of differential geometry, the geodesic in $D$%
-dimensional Riemannian space is considered as {\sl special kind of a curve,
having the following properties}.

\noindent (i) {\sl Extremality}. The distance $(2\sigma )^{1/2}$, measured
along the geodesic between two points is the shortest (extremal) as compared
with the distance measured along other curves.

\noindent (ii) {\sl Definiteness}. Any two different points of the geodesic
determine uniquelly the geodesic, passing through these points.

\noindent (iii) {\sl Minimality of the section} (one-dimensionality). Any
section of the geodesic consists of one point.

At the conventional approach the property (ii) is a corollary of the
property (i) (for rather small regions of the space), but the property (iii)
is the property of any curve (but not only of geodesic).

In T-geometry the geodesic is considered as a special kind of the tube,
degenerating into a line. Then the properties (ii) and (iii) are supposed to
be fulfilled. The property (i) is not defined, because the concept of line
is not defined.

Let us try to determine the geodesic as the tube, having the properties of
definiteness and of the section minimality at the same time.

\begin{opred}
\label{d3.1.10}. The tube ${\cal T}({\cal P}^n)$ has the definiteness
property, if for any $(n+1)$-point tube $\sigma $-basis ${\cal Q}^n\subset
{\cal T}({\cal P}^n)$ (or for any multivector $\overrightarrow{{\cal Q}^n}%
\subset {\cal T}({\cal P}^n)$ of unvanishing length) the following condition
is fulfilled
\begin{equation}
{\cal T}({\cal Q}^n)={\cal T}({\cal P}^n)  \label{a2.49}
\end{equation}
\end{opred}

\begin{opred}
\label{d3.1.11}. The tube ${\cal T}({\cal P}^n)$ has the minimality section
property, if $\forall P\in {\cal T}({\cal P}^n)$
\begin{equation}
{\cal S}_{n;P}({\cal T}({\cal P}^n))=\{P\},\qquad \forall P\in {\cal P}^n
\label{a2.50}
\end{equation}
\end{opred}

\begin{opred}
\label{d3.1.12}. $\sigma $-space is extremal on the tube ${\cal T}({\cal P}%
^n)$, if for ${\cal T}({\cal P}^n)$ the conditions of definiteness and
section minimality are fulfilled.
\end{opred}

\begin{opred}
\label{d3.1.13}. $\sigma $-space is extremal on the set ${\cal T}$ of tubes $%
{\cal T}({\cal P}^n)$, if it is extremal on any tube of the set ${\cal T}$.
\end{opred}

\begin{opred}
\label{d3.1.14}. $\sigma $-space is extremal in the $n$th order, if it is
extremal on all $n$th order tubes ${\cal T}({\cal P}^n)$.
\end{opred}

\begin{opred}
\label{d3.1.15}. The tube ${\cal T}({\cal P}^n)$ is the geodesic tube ${\cal %
L}({\cal P}^n)$, if the $\sigma $-space is extremal on the tube ${\cal T}(%
{\cal P}^n)$.
\end{opred}

\section{Euclidean space as a special case of $\sigma $-space.}

\begin{opred}
\label{d3.2.1}. $n$-dimensional Euclidean space $E_n$ is a set ${\Bbb R}^n$
of all ordered sets $x=\{x_1,x_2,\ldots x_n\}$ of $n$ real numbers on which
for $\forall x\in {\Bbb R}^n$, $\forall y\in {\Bbb R}^n$ is given the real
function $\sigma $:
\begin{equation}
\sigma (x,y)={\frac 12}\sum_{i,k=1}^ng^{ik}(x_i-y_i)(x_k-y_k),\qquad g^{ik}=%
\hbox{const},\qquad i,k=1,2,\ldots n  \label{a3.1}
\end{equation}
\begin{equation}
\det \parallel g_{ik}\parallel =(\det \parallel g^{ik}\parallel )^{-1}\neq 0
\label{a3.2}
\end{equation}
\end{opred}

The function $\sigma $ is called the world function or simply $\sigma $%
-function. $n$-dimensional Euclidean space $E_n$ is at the same time a $%
\sigma$-space $E_n=\{\sigma , {\Bbb R}^n\}$.

{\it Remark.} The given definition is equivalent to the definition of $n$%
-dimensional Euclidean space $E_{n}$ as $n$-dimensional linear space ${\Bbb R%
}^n$ of vectors $x=\{x_1,x_2,\ldots ,x_n\}\in {\Bbb R}^n$ with given on it
the scalar product $(x.y)$ of vectors $x,y\in {\Bbb R}^n$
\begin{equation}
(x.y)=\sum^{n}_{i,k=1}g^{ik}x_{i}y_{k}=\sigma(0,x)+\sigma(0,y)-\sigma(x,y),
\qquad g^{ik}=\hbox{const},\qquad i,k=1,2,\ldots n  \label{a3.3}
\end{equation}
where $\sigma$ is given by the relation (\ref{a3.1}).

The Euclidean space $E_n=\{\sigma ,\Omega \}$, $({\Omega ={\Bbb R}}^n)$,
considered as $\sigma$-space, have the following properties
\begin{equation}
\exists {\cal P}^n\subset\Omega,\qquad F_n({\cal P}^n)\ne 0,\qquad
F_{n+1}(\Omega ^{n+2})=0,  \label{a3.4}
\end{equation}
\[
\sigma (P,Q)={\frac{1}{2}}\sum^{n}_{i,k=1}g^{ik}({\cal P}^{n}) [\Gamma
(P_{0},P_{i},P)-\Gamma (P_{0},P_{i},Q)]
\]
\begin{equation}
\times [\Gamma (P_{0},P_{k},P)-\Gamma (P_{0},P_{k},Q)], \qquad \forall
P,Q\in\Omega  \label{a3.5}
\end{equation}
\begin{equation}
\Gamma (P_{0},P,Q)=\sum^{n}_{i,k=1}g^{ik}({\cal P}^{n})\Gamma
(P_{0},P_{i},P) \Gamma (P_{0},P_{k},Q),\qquad \forall P,Q\in\Omega ,
\label{a3.6}
\end{equation}
where ${\cal P}^n$ is some $(n+1)$-point tube $\sigma$-basis in $\Omega =%
{\Bbb R}^n$ in the sense of the definition (\ref{a2.32}), i.e. $F_n({\cal P}%
^n)\ne 0$, and the quantities $\Gamma (P_{0},P_{k},P)$ are defined by the
relations (\ref{a2.3}). $(n+1)$-point tube $\sigma$-basis ${\cal P}^n$
corresponds to the basis of $n$ vectors
\begin{equation}
{\bf e}_i={\bf P}_{0}{\bf P}_{i},\qquad P_i\in{\cal P}^n, \qquad i
=1,2,\ldots n  \label{a3.7}
\end{equation}
and
\begin{equation}
x_{i}= x_{i}(P) = ({\bf P}_{0}{\bf P} .{\bf e}_{i})=\Gamma (P_{0},P,P_{i}),
\qquad i = 1,2,\ldots n,\qquad\forall P\in\Omega  \label{a3.8}
\end{equation}
are covariant coordinates of the vector ${\bf P}_{0}{\bf P}$ in this basis.
The quantities
\begin{equation}
g_{ik}= g_{ik}({\cal P}^{n}) = ({\bf e}_{i}.{\bf e}_{k})= \Gamma
(P_{0},P_{i},P_{k}),\qquad i,k=1,2,\ldots n  \label{a3.9}
\end{equation}
are covariant components of the metric tensor in this basis ${\cal P}^n$. As
far as ${\cal P}^n$ is the tube $\sigma$-basis, then according to (\ref
{a2.32}), (\ref{a2.33}) the following condition is fulfilled
\begin{equation}
F_n({\cal P}^n)\equiv \det \parallel g_{ik}({\cal P}^n)\parallel \neq 0,
\qquad i,k=1,2,\ldots n  \label{a3.10}
\end{equation}
and one can determine the contravariant components $g^{ik}=g^{ik} ({\cal P}%
^n)$ of the metric tensor by means of the relation
\begin{equation}
\sum^{n}_{k=1}g_{ik}({\cal P}^{n})g^{kl}({\cal P}^{n}) = \delta ^{l}_{i},
\qquad i,l=1,2,\ldots n  \label{a3.11}
\end{equation}

Conditions (\ref{a3.5}) and (\ref{a3.6}) are equivalent, as it follows from (%
\ref{a2.3}).

\begin{opred}
\label{d3.2.4}. $\sigma $-space $V=\{\sigma ,\Omega \}$ have the structure
of $n$-dimensional Euclidean space on $\Omega $, if there exists such a $%
(n+1)$-point tube $\sigma $-basis ${\cal P}^n\subset \Omega $, that $\forall
P,$ $\forall Q\in \Omega $ the condition (\ref{a3.6}) is fulfilled.
\end{opred}

$\sigma $-space $V$, having the structure of $n$-dimensional Eucliden space
may be not Euclidean, because one-to-one correspondence between the points $%
P\in \Omega $ and their coordinates $x\in {\Bbb R}^{n}$ may not exist. For
instance, two different points $P$ and $P^\prime$ may have similar
coordinates and be mapped on one point $x$ of the Euclidean space $%
E_n=\{\sigma ,{\Bbb R}^n\}$.

Finally, the third property of the Euclidean space $E_n=\{\sigma , {\Bbb R}%
^n\}$ is formulated as follows. The relation
\begin{equation}
\Gamma (P_0,P_i,P)=x_i, \qquad x_i\in {\Bbb R},\qquad i=1,2,\ldots n,
\label{a3.12}
\end{equation}
considered as equations for determination of $P\in\Omega ={\Bbb R}^n$,
always have one and only one solution.

Let us note that all three conditions are written in $\sigma $-immanent
form. They are necessary properties of the Euclidean space. In this
connection one can put the question whether these conditions are also
sufficient conditions for the $\sigma $-space $\{\sigma ,\Omega \}$ were the
Euclidean space. The following theorem answers this question.

\begin{theorem}
For the $\sigma $-space $\{\sigma ,\Omega \}$ were $n$-dimensional Euclidean
space, It is necessary and sufficient that the conditions (\ref{a3.4}), (\ref
{a3.5}) and (\ref{a3.12}) be fulfilled. \label{c2}
\end{theorem}

{\it Proof.} {\it Necessity} of conditions (\ref{a3.4}), (\ref{a3.5}) and (%
\ref{a3.12}) is tested by the direct substitution of world function $\sigma$
for $n$-dimensional Euclidean space $E_n=\{\sigma , {\Bbb R}^n\}$.

{\it Sufficiency}. Let ${\cal P}^n$ be some $(n+1)$-point tube $\sigma$%
-basis in $\sigma $-space $\{\sigma ,\Omega \}$ and $P,Q\in\Omega$ be two
arbitrary points. Let us introduce their covariant coordinates in ${\cal P}%
^n $ by means of the relations of type (\ref{a3.8})
\begin{equation}
x_{i}= \Gamma (P_{0},P_{i},P),\qquad y_{i}= \Gamma (P_{0},P_{i},Q), \qquad i
= 1,2,\ldots n,\qquad\forall P,Q\in\Omega  \label{a3.13}
\end{equation}
Then the relation (\ref{a3.6}) is rewritten in the form
\begin{equation}
(x.y)=\sum^{n}_{i,k=1}g^{ik}x_{i}y_{k}, \qquad g^{ik}=g^{ik}({\cal P}^{n})=%
\hbox{const},\qquad i,k=1,2,\ldots n  \label{a3.14}
\end{equation}
In virtue of the condition (\ref{a3.12}) any point $P\in\Omega$ corresponds
to one and only one point $x\in {\Bbb R}^n$ and vice versa. In other words,
the $\sigma $-space $\{\sigma ,\Omega \}$ is isometric to $n$-dimensional
Euclidean space $E_n=\{\sigma , {\Bbb R}^n\}$.

{\it Corollary of the theorem.} $n$-dimensional Euclidean space and all its
properties can be described $\sigma$-immanently (i.e. in terms of the world
function). In other words, the T-geomtry is rich and pithy enough to contain
Euclidean geometry as a special case, when the world function is restricted
by $\sigma$-immanent relation (\ref{a3.5}), or by equivalent relation (\ref
{a3.6}). The Riemanian geometry can be presented in the $\sigma$-immanent
form \cite{R90} also. This may be interpreted in the sense that T-geometry
contains the Riemannian geometry as a special case. T-geometry is
informative enough to contain other geometries. Apparently, it is rather
difficult to construct a geometry which would not be contained in
T-geometry. The fact is that that practically any geometry may be considered
as a result of a deformation (variation of the world function) of $\sigma$%
-subspace of the Euclidean space of rather high dimensionality. At such a
variation the pithiness of geometry does not reduce, because the number of
tubes does not reduce. This number may only increase under deformation,
because any tube of the Euclidean space, determined by many $\sigma$-bases,
is splitted, in general, to several different tubes.

Pithiness of geometry (i.e. the number of geometric objects, suppositions
and theorems) depends not only on axioms of the geometry, it depends also on
development of the mathematical technique of the geometry. Conventionally
the metric geometry is considered as the geometry which is less pithy as
compared with the Euclidean geometry. One connects usually the pithiness of
the metric geometry with constraints (\ref{a1.3}), (\ref{a1.4}), which are
essential for construction of the shortest and geometric objects, connected
with it. Essentially, the pithiness of the metric geometry is connected with
its mathematical technique. Using more effective mathematical technique,
connected with the classification (\ref{a1.4}), the pithiness of metric
geometry increases even under removing constraints (\ref{a1.3}), (\ref{a1.4}%
) on metric.

T-geometry has such a dignity as insensitivity of its mathematical technique
to continuity of the set, where the geometry is given. For instance, let us
take 100 points ${\cal P}^{99}$ of the three-dimensional Euclidean space and
try to study geometry of this set of points. Using usual way, one should
consider Euclidean space on the set ${\Bbb R}^3$, introduce a coordinate
system, remove all points except for ${\cal P}^{99}$ and begin to study the
way of embedding the set ${\cal P}^{99}$ in the Euclidean space, starting
from coordinates of its points. From point of view of T-geometry one should
study the set ${\cal P}^{99}$, imposing the constraint (\ref{a3.5}) on
metric and removing (\ref{a3.12}). Thus, approach of T-geometry appears to
be local in the sense that the geometry of the set ${\cal P}^{99}$ is
studied, but not the way of embedding the set in the Euclidean space.

Giving up of constraints (\ref{a3.12}) leads to a violation of the mapping $%
\Omega \to {\Bbb R}^n $ reversibility. In particular, it is possible such a
case, when the $\sigma $-space $\{\sigma ,\Omega \}$ appears to be a $\sigma
$-subspace of the Euclidean space $E_n$.

\begin{opred}
\label{d3.2.2}. The Euclidean $\sigma $-space $E^{\prime }=\{\sigma ,\Omega
^{\prime }\}$ is the $\sigma $-space which can be isometrically embedded in
the Euclidean space. $n$-dimensional Euclidean $\sigma $-space $E_n^{\prime
}=\{\sigma ,\Omega ^{\prime }\}$ is $\sigma $-space which can be
isometrically embedded in $n$-dimensional Euclidean space $E_n=\{\sigma ,%
{\Bbb R}^n\}$, but cannot be isometrically embedded in $(n-1)$-dimensional
Euclidean space $E_{n-1}=\{\sigma ,{\Bbb R}^{n-1}\}$.
\end{opred}

$n$-dimensional Euclidean $\sigma$-space is a $\sigma$-subspace of $n$%
-dimensional Euclidean space $E_n=\{\sigma , {\Bbb R}^n\}$.

From the tube definition (\ref{a2.35}) and the condition (\ref{a3.4}) it
follows that $n$-dimensio\-nal Euclidean $\sigma$-space is the $n$th order
tube ${\cal T}({\cal P}^n)=\Omega$, generated by any $(n+1)$-point tube $%
\sigma$-basis ${\cal P}^n\subset\Omega$, the condition of the tube section
minimality (\ref{a2.50}) being fulfilled. Then the following theorem takes
place.

\begin{theorem}
\label{t3.2.2}. Let ${\cal P}^n$ be $(n+1)$-point tube $\sigma $-basis in
the $\sigma $-space $V\{\sigma ,\Omega \}$. For the tube ${\cal T}({\cal P}%
^n)$ be $n$-dimensional Euclidean $\sigma $-space it is necessary and
sufficient, that

\noindent (1)\quad $\sigma $-space ${\cal T}({\cal P}^n)$ have the structure
of $n$-dimensional Euclidean space on ${\cal T}({\cal P}^n)$,

\noindent (2)\quad Section of ${\cal T}({\cal P}^n)$ be minimal at any
point:
\[
{\cal S}_{n;P}({\cal T}({\cal P}^n))=\{P\},\qquad \forall P\in {\cal T}(%
{\cal P}^n).
\]
\end{theorem}

\section{Triangle axiom as a condition of the first order tube degeneration}

Let us study constraints, imposed on $\sigma$-space by the triangle
inequality (\ref{a1.4}). Let us consider segment ${\cal T}_{[P_0P_1]}$ of
the tube ${\cal T}_{P_0P_1}$, contained between basic points $P_0,P_1$. It
is described by equations (\ref{c.1.23}), (\ref{a2.41}).

For continuous $\sigma$-space the tube segment ${\cal T}_{[P_0P_1]}$ is some
surface, containing points $P_0,P_1$. This surface ${\cal T}_{[P_0P_1]}$ and
the region outside the surface are described by the equation
\begin{equation}
S_2(P_0,R,P_1)\equiv \rho (P_0,R)+\rho (R,P_1 )- \rho (P_0,P_1 )\ge 0,
\label{a4.1}
\end{equation}
where $R$ is the running point. Thus, the triangle inequality is fulfilled
on the surface ${\cal T}_{[P_0P_1]}$ and outside it. The region inside the
surface ${\cal T}_{[P_0P_1]}$ associates with the inequality $%
S_2(P_0,R,P_1)\le 0$, that corresponds to a violation of the triangle axiom.
In other words, in the metric space the first order tube segment ${\cal T}%
_{[P_0P_1]}$ has no inner points. It means degeneration of the tube into a
line, or into a surface which has no inner points. In this sense the metric
geometry (i.e. geometry generated by the metric space) is degenerated
geometry.

In the case, when all first order tubes ${\cal T}_{P_0P_1}$ degenerate into
corresponding basic points $P_0,P_1$, the triangle inequality (\ref{a1.4})
takes the form of a strong inequality
\begin{equation}
\rho (P_0,R)+\rho (R,P_1)> \rho (P_0,P_1),\qquad P_0\ne R\ne P_1\ne P_0,
\qquad \forall P_0,P_1,R\in \Omega .  \label{a4.2}
\end{equation}
In this case it seems to be reasonable to call the T-geometry
ultradegenerated.

{\it Example} 4.1. Let us consider two different $\sigma$-spaces (and two
T-geometries) on the unit sphere.
\begin{equation}
\Omega =\left\{ {\bf x}\left| |{\bf x}|^2\leq 1\right. \right\} \subset
{\Bbb R}^3,\qquad {\bf x}=\left\{ x^1,x^2,x^3\right\} \in {\Bbb R}^3,\qquad |%
{\bf x|}^2\equiv \sum\limits_{i=1}^3\left( x^i\right) ^2  \label{a4.3}
\end{equation}
$\sigma$-space $V_E=\left\{ \sigma _E,\Omega \right\} $ generates the proper
Euclidean geometry
\begin{equation}
\sigma _E:\quad \Omega \times \Omega \rightarrow [0,\infty ) \subset {\Bbb R}%
,\qquad \sigma _E({\bf x},{\bf x}^{\prime })=\frac 12|{\bf x}-{\bf x}%
^{\prime }|^2,\qquad {\bf x},{\bf x}^{\prime }\in \Omega ,  \label{a4.4}
\end{equation}
$\sigma$-space $V=\left\{ \sigma ,\Omega \right\} $ generates T-geometry on
the same set $\Omega $ by means of relations
\begin{equation}
\sigma :\quad \Omega \times \Omega \rightarrow [0,\infty )\subset {\Bbb R},
\qquad \sigma ({\bf x},{\bf x}^{\prime })=2\left( \arcsin \sqrt{\frac{\sigma
_E({\bf x},{\bf x}^{\prime })}2}\right) ^2,\qquad {\bf x},{\bf x}^{\prime
}\in \Omega  \label{a4.5}
\end{equation}
Along with the two $\sigma$-spaces in the sphere $\Omega $ one considers
their $\sigma$-subspaces $V_{Es}=\left\{ \sigma _E,\Sigma \right\} $ and $%
V_s=\left\{ \sigma ,\Sigma \right\} $ on the sphere surface $\Sigma
=\partial \Omega =\left\{ {\bf x}\left| |{\bf x}|^2=1\right. \right\}
\subset \Omega $. As far as $\Sigma $ is a subset of the set $\Omega $, the
'-geometries $V_{Es}=\left\{ \sigma _E,\Sigma \right\} $ and $V_s=\left\{
\sigma ,\Sigma \right\} $ are generated by T-geometries $T_E$ ¨ $T$.

Let us design the tubes in $\sigma $-space $V_E$ by means of the symbol $%
{\cal L}$, the tubes in the $\sigma $-space $V$ are denoted by the symbol $%
{\cal T}$. The first order tubes ${\cal L}_{AB}\subset \Omega ,\;\;\left(
A,B\in \Sigma \right) $ are straight lines in $\Omega $, and they are formed
by two basic points $A,B$ in $\Sigma $. In other words, T-geometry $V_E$ is
degenerated in the first order in $\Omega $, and it is ultradegenerated in
the first order in $\Sigma $. The first order tubes ${\cal T}_{AB}\subset
\Omega ,\;\;\left( A,B\in \Sigma \right) $ are nondegenerated tubes in $%
\Omega .$ They are surfaces, formed by a rotation of unit radius circles,
passing through points $A,B\in \Sigma $, around the axis ${\cal L}_{AB}$.
(see. Figure 4.1). These tubes tangent the sphere $\Sigma $ along the
circles of maximal radius. The segment ${\cal T}_{[AB]}$ of the tube between
the points $A,B\in \Omega $ is found inside the sphere $\Omega ,$ whereas
the remaining part of the tube ${\cal T}_{AB}$ is found outside the inner
part $\Omega \backslash \Sigma $ of the sphere. As a result the segment $%
{\cal T}_{[AB]}$ of the tube ${\cal T}_{AB}\subset \Sigma ,\;\;\left( A,B\in
\Sigma \right) $ is the shortest in the $\sigma $-space $V_s=\left\{ \sigma
,\Sigma \right\} $. This shortest on the sphere surface $\Sigma $ connects
points $A,B\in \Sigma $. The remaining part of the tube ${\cal T}%
_{AB}\subset \Omega $ is a continuation of the segment ${\cal T}%
_{[AB]}\subset \Sigma $. In other words, $\sigma $-space $V_s=\left\{ \sigma
,\Sigma \right\} $ generates the degenerated in the first order T-geometry
on $\Sigma .$ Thus, T-geometry $V=\left\{ \sigma ,\Omega \right\} $ is
nondegenerated in $\Omega $ and it is degenerated in $\Sigma $. The second
order tube ${\cal T}_{ABC}\subset \Sigma $ consists of three points $%
A,B,C\subset \Sigma $, and T-geometry in $V_s=\left\{ \sigma ,\Sigma
\right\} $ is ultradegenerated in the second order geometry in $\Omega $.

On the other hand, T-geometry in $V_s=\left\{ \sigma ,\Sigma \right\} $ on
the sphere surface $\Sigma $ can be constructed on the basis of Euclidedan
geometry in $V_E=\left\{ \sigma _E,\Omega \right\} $. To construct $%
V_s=\left\{ \sigma ,\Sigma \right\} $ on the basis of $V_E=\left\{ \sigma
_E,\Omega \right\} $, one can use extremal properties of geodesics.

Let us consider the second order tube ${\cal L}_{ABC}\subset \Omega $, $%
(A,B,C\in\Omega )$. This tube is a two-dimensional plane, passing through
the points $A,B,C\in \Omega $. In $V_{Es}=\left\{ \sigma _E,\Sigma \right\} $
the second order tube ${\cal L}_{ABC}\subset \Sigma $ has the form of a
circle, passing through points $A,B,C\in\Sigma $.

To construct internal geometry in $V_s=\left\{ \sigma ,\Sigma \right\} $, it
is necessary to determine the metric $\rho\left( A,B\right) = \sqrt{2\sigma
\left( A,B\right) },\;\;A,B\in \Sigma $ on $\Sigma \times\Sigma $, using the
following way.
\begin{equation}
\rho\left( A,B\right) =\inf\limits_{C\in\Sigma ,C\neq A,C\neq B}l_C\left(
A,B\right) ,\qquad A,B\in\Sigma ,  \label{a4.6}
\end{equation}
where $l_C\left( A,B\right)\subset [0,\infty ) $ is the length of the curve $%
{\cal L}_{ABC}\subset\Sigma $ between the ponts $A,B$. Let us order the
points $R\in {\cal L}_{ABC}\subset\Sigma $ of the curve ${\cal L}%
_{ABC}\subset\Sigma $, solving the equation
\begin{equation}
\rho _E\left( A,R\right) \left( 1-\tau \right) =\rho _E\left( B,R\right)
\tau ,\qquad R\in {\cal L}_{ABC}\subset \Omega ,\qquad \tau \in {\Bbb R}
\label{a4.7}
\end{equation}
Solution of this equation determines $R=R_{AB}\left( \tau ,C\right) \in
{\cal L}_{ABC}\subset\Sigma $ as a function of the parameter $\tau \in [0,1]$
and of the point $C\in\Sigma $. Therewith $A=R_{AB}\left(0,C\right)$, $%
B=R_{AB}\left(1,C\right)$. The function $R=R_{AB}\left( \tau ,C\right)$ has
two branches $R=R_1\left( \tau ,C\right) $ and $R=R_2\left( \tau ,C\right)$.
It should change the branch with the less values of $\rho _E(A,R_{AB}\left(
\tau ,C\right))$ and determine $l_C\left( A,B\right) $ by means of the
relation
\begin{equation}
l_C\left( A,B\right) =\int\limits_0^1\left[ \frac{d\rho _E\left(R_{AB}\left(
\tau ,C\right) ,R_{AB}\left( \tau ^{\prime },C\right) \right) }{d\tau
^{\prime }}\right] _{\tau ^{\prime }=\tau }d\tau  \label{a4.8}
\end{equation}

Substituting (\ref{a4.8}) in (\ref{a4.6}), one obtains
\begin{equation}
\rho \left( A,B\right) =2\arcsin \frac{S_E\left( A,B\right) }2,\qquad A,B\in
\Sigma  \label{a4.9}
\end{equation}
that corresponds to (\ref{a4.5}).

This example shows that, using Euclidean geometry inside $\Omega $, one can
consrtuct internal metric (and T-geometry) on the surface $\Sigma $ of the
sphere $\Omega $. The extremal properties of geodesics (shortests) and the
second order tubes ${\cal L}_{ABC}\subset \Omega $ are used essentially for
construction of geometry on $\Sigma $. The second order tubes ${\cal L}%
_{ABC}\subset \Omega $ generate a system of curves on the sphere surface $%
\Sigma $. One chooses those among them which have the minimal length between
points $A,B$.

Apparently, with proper stipulations this procedure of constructing metric
space on the sphere surface with Euclidean space inside sphere can be
generalized on the case of arbitrary body and arbitrary metric space inside
it.

Thus, classification of finite metric spaces, using the series of mappings (%
\ref{a1.5}), appears to be a very effective method of studying the metric
space. This method admits to describe the Euclidean geometry and the
Riemannian one in terms of only metric (world function). Geometries,
constructed on the basis of this classification do not use concept of
continuity. They are insensitive to discreteness or continuity of the space.
Concept of continuity may be introduced on the basis of metric (world
function) by means of a proper parametrization of extremal tubes \cite{R90}.



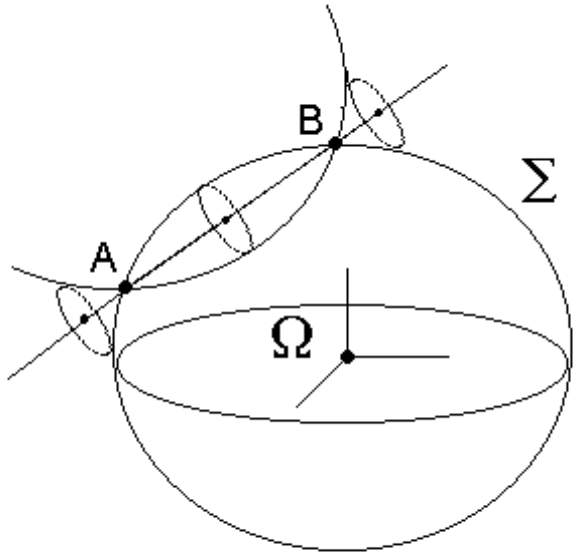
\end{document}